\documentclass[11pt,draft]{article}
\usepackage{amsfonts}
\usepackage{mathrsfs}
\usepackage{amsmath,amsfonts,mathrsfs,amssymb,color}
\usepackage{indentfirst}

\numberwithin{equation}{section}

\newtheorem{theo.}{\quad\, Theorem}[section]
\newtheorem{defi.}{\quad\, Definition}[section]
\newtheorem{lemm.}{\quad\, Lemma}[section]
\newtheorem{coro.}{\quad\, Corollary}[section]

\begin{document}

\title{
Spectrum structure for eigenvalue problems involving mean curvature
operators in Euclidean and Minkowski spaces
$^*$ }
\author{
Ruyun Ma$^{a}$ \ \ \ \
Hongliang Gao$^{b}$ \ \ \ \ Tianlan Chen$^c$
\\
 {\small $^{a,b,c}$Department of Mathematics, Northwest
Normal University, Lanzhou 730070, P R China}\\
}
\date{} \maketitle
\footnote[0]{E-mail address: mary@nwnu.edu.cn (R. Ma), \ \ gaohongliang101@163.com (H.Gao), \ \
chentianlan511@126.com (T. Chen).
  \
} \footnote[0] {$^*$Supported by the  NSFC (No.11361054, No.11201378), SRFDP(No.20126203110004), Gansu provincial National
Science Foundation of China (No.1208RJZA258). }
\begin{abstract}
In this paper, we are concerned with quasilinear Dirichlet problem
$$
\left\{
\aligned
&-\Big(\frac{u'(x)}{\sqrt{1+\kappa (u'(x))^2}}\Big)'=\lambda u(x), \ \
 \ \ \ 0<x<1,\\
&u(0)= u(1)=0,\\
\endaligned
\right.
\eqno (P)
$$
where $\kappa\in (-\infty, 0)\cup (0, \infty)$ is a constant.  We show that any nontrivial solution $ u$ of (P) has only finite many of simple zeros in $[0,1]$, all of humps of $u$ are same, and the first hump is symmetric around the middle point of its domain. We also describe the global structure of the set of nontrivial solutions of (P).

\end{abstract}


{\small\bf MR(2000)\ \ \ 35K59, \ 35K93}

\baselineskip 24pt

\section{Introduction}

This paper studies the {\it spectrum structure} of the quasilinear Dirichlet problem
$$
\left\{
\aligned
&-\Big(\frac{u'(x)}{\sqrt{1+\kappa (u'(x))^2}}\Big)'=\lambda u(x), \ \
 \ \ \ 0<x<1,\\
&u(0)= u(1)=0,\\
\endaligned
\right.
\eqno (1.1)
$$
where
$$'=D:=\frac{d}{dx},
$$
and $\kappa\neq 0$ is a constant. This problem establishes a quasilinear continuum deformation between the linear eigenvalue problem
$$
\left\{
\aligned
&-u''(x)=\lambda u(x), \ \
 \ \ \ 0<x<1,\\
&u(0)= u(1)=0,\\
\endaligned
\right.
\eqno (1.2)
$$
and (1.1), which is a quasilinear problem associated to
the one-dimensional Euclidean-curvature operator as $\kappa>0$ and to the one-dimensional Minkowski-curvature operator as $\kappa<0$, respectively.

The quasilinear problem (1.1) can be equivalently written as
$$
\left\{
\aligned
&-\frac{u''(x)}{[1+\kappa (u'(x))^2]^{3/2}}=\lambda u(x), \ \
 \ \  \ 0<x<1,\\
&u(0)= u(1)=0,\\
\endaligned
\right.
\eqno (1.3)
$$
or
$$
\left\{
\aligned
&-u''(x)=\lambda u [1+\kappa (u'(x))^2]^{3/2}, \ \
 \ \  \ 0<x<1,\\
&u(0)= u(1)=0.\\
\endaligned
\right.
\eqno (1.4)
$$

Cano-Casanova,  L\'{o}pez-G\'{o}mez and Takimoto [1] used the
standard techniques from bifurcation theory [2-5] to study  the
existence of positive solutions for (1.4) with $\kappa>0$. Their
main results can be summarized in the following list:

$\bullet$ Problem (1.4) has a positive solution if and only if
$$8B^2<\lambda<\pi^2, \ \ \ \ \ \ B=\int^1_0 \frac{d\theta}{\sqrt{\theta^{-4}-1}}.
\eqno (1.5)
$$

$\bullet$ The positive solution of (1.4) is unique if it exists.
Subsequently, we denote it by $u_\lambda$.

$\bullet$ $u_\lambda$ is symmetric around $1/2$ for all $\lambda$ satisfying (1.5).

$\bullet$ $u_\lambda'(0)$ is strictly decreasing in $(8B^2, \pi^2)$.

$\bullet$ $u_\lambda$ satisfies
$$
    \lim_{\lambda\downarrow 8B^2}u'_\lambda(0)=\infty, \ \ \ \ \ \
    \lim_{\lambda\downarrow 8B^2}||u_\lambda||_\infty=\frac 1{2B\sqrt{\kappa}},
\eqno (1.6)
$$
and
$$\lim_{\lambda\uparrow \pi^2}||u_\lambda||_\infty=\lim_{\lambda\uparrow \pi^2}u_\lambda(1/2)=0.
\eqno (1.7)
$$

$\bullet$ The point-wise limit $u_{8B^2}:=\lim_{\lambda\downarrow 8B^2} u_\lambda$ satisfies
$$u_{8B^2}'(0)=-u_{8B^2}'(1)=\infty.
\eqno (1.8)
$$
Their results are motivated by the pioneering results in Nakao [6] and have been extended to the more general case
$$
\left\{
\aligned
&-\Big(\frac{u'(x)}{\sqrt{1+\kappa (u'(x))^2}}\Big)'=\lambda V(x)u(x), \ \
 \ \ \ 0<x<1,\\
&u(0)= u(1)=0,\\
\endaligned
\right.
$$
by Cano-Casanova,  L\'{o}pez-G\'{o}mez and Takimoto [7], where the
the weight function $V(x)\geq 0$ and $V(x)\not\equiv 0$ in $[0,1]$.

\vskip 3mm

It is well-known that (1.2) has a sequence of eigenvalues
$$\lambda_n=n^2\pi^2, \ \ \ \ n\in \{1,\, 2,\, \cdots\}.
$$
For each $n\in \mathbb{N}$, $\lambda_n$ is simple, and its eigenfunction $\varphi_n=\sin n\pi x$ has exactly $n-1$ simple zeros in $(0,1)$.

\vskip 3mm

  Of course, the natural question is what would happen if we consider the sign-changing solutions of (1.1) under $\kappa\neq 0$?
  It is the purpose of this paper to study the {\it higher eigenvalue} case for (1.1).
   More precisely, we shall show the following:
\vskip 3mm

\noindent{\bf Theorem 1.1} Assume that $\kappa>0$. Let
$$S^\nu_n:=\{u\in C^1_0[0,1]\,|\,u \, \text{has exactly}\; n-1\;
             \text{simple zeros in}\ (0,1), \ \text{and}\ \nu  u\
             \text{is positive near}\ 0\}.
             $$
for $n\in N$ and $\nu\in \{+,-\}$.  Then

1. Problem (1.4) has a $S^\nu_n$-solution  if and only if
$$8n^2B^2<\lambda<n^2\pi^2, \ \ \ \ \ \ B=\int^1_0 \frac{d\theta}{\sqrt{\theta^{-4}-1}}.
\eqno (1.9)
$$

2. The $S^\nu_n$-solution of (1.4) is unique if it exists.
Subsequently, we denote it by $u_\lambda$.

3.  All of the humps of $u_\lambda$ are same, and the first hump is symmetric around $\frac 1{2n}$ for all $\lambda$ satisfying (1.9).

4. $u_\lambda'(0)$ is strictly decreasing in $(8n^2B^2, n^2\pi^2)$.

5. $u_\lambda$ satisfies
    $$
    \lim_{\lambda\downarrow 8n^2B^2}||u_\lambda||_\infty=\frac 1{2nB\sqrt{\kappa}},
$$
and
$$\lim_{\lambda\uparrow n^2\pi^2}||u_\lambda||_\infty=\lim_{\lambda\uparrow n^2\pi^2}u_\lambda(\frac 1{2n})=0.
$$

6.  The point-wise limit
$$u_{8n^2B^2}:=\lim_{\lambda\downarrow 8n^2B^2} u_\lambda$$
provides us with a solution of (1.4) at $\lambda=8n^2B^2$ which is singular at $t=\frac jn$, $j=0, 1, \cdots n$,  in the sense that
$$(-1)^j u_{8B^2}'(\frac jn)=\infty, \ \ \ \ \ j=0, 1, \cdots n.
\eqno (1.10)$$

\vskip 3mm

\noindent{\bf Theorem 1.2} Assume that $\kappa<0$. Then

(1) For $n\in N$ and $\nu\in \{+,-\}$, (1.4) has a $S^\nu_n$-solution  if and only if
$$n^2\pi^2<\lambda<\infty.
\eqno (1.11)
$$

(2) The $S^\nu_n$-solution of (1.4) is unique if it exists.
Subsequently, we denote it by $u_\lambda$.

(3) All of the humps of $u_\lambda$ are same, and the first hump is symmetric around $\frac1{2n}$ for all $\lambda$ satisfying (1.11).

(4) $u_\lambda'(0)$ is strictly increasing in $(n^2\pi^2, \infty)$.

(5) $u_\lambda$ satisfies
    $$
    \lim_{\lambda\uparrow \infty}||u_\lambda||_\infty\leq\frac 1{2n\sqrt{-\kappa}},
$$
and
$$\lim_{\lambda\downarrow n^2\pi^2}||u_\lambda||_\infty=\lim_{\lambda\downarrow n^2\pi^2}u_\lambda(\frac 1{2n})=0.
$$

(6) The point-wise limit
$$u_{\infty}:=\lim_{\lambda\uparrow \infty} u_\lambda$$
satisfies
$$
u'_{\infty}(\frac jn):=(-1)^j\frac 1{\sqrt{-\kappa}}, \ \ \ \ j\in\{0,\, 1,\, \cdots, n\}.
$$

   \vskip 4mm

  Throughout this paper we use the following notations and conventions.
Given a function $u\in C[a, b]$ it is said that $u \succ 0$ if
$u\geq 0$ but $u\not\equiv 0$.
For any $V\in C[a, b]$, we shall denote by $\sigma[-D^2 + V;(a,b)]$
the principal eigenvalue of $-D^2 + V$ under homogeneous Dirichlet
boundary conditions in the interval $(a,b)$. In particular,
$\sigma[-D^2, (0,1)]=\pi^2$.

Let $E:=C^1_0[0,1]$ with the normal
$$||u||_1:=\max\{||u||_\infty, ||u'||_\infty\}.
$$
 Then $S^\nu_n$ is open in $E$ for each $n\in N$ and $\nu\in\{+, -\}$.

\vskip 3mm

     For the related results on the quasilinear problem (1.1) and its more general case (including the  higher dimensional case), see D. Gilbarg and N. S. Trudinger [8], S.-Y. Cheng and S.-T. Yau [9],
R. Bartnik and L. Simon [10],   C. Bereanu, P. Jebelean and J. Mawhin [11],  C. Bereanu, P. Jebelean, P. J. Torres [12-13], M. F. Bidaut-V\'{e}ron and A. Ratto [14],
I. Coelho, C. Corsato, F. Obersnel and P. Omari [15],  R. L\'{o}pez [16], J. Mawhin [17], A. E. Treibergs [18], A. Azzollini [19],  H. Pan and R. Xing [20] and references therein.

\vskip 3mm
     The contains of this paper have been distributed as follows. In Section 2, we shows that any nontrivial solutions of (1.4) belongs to $S^\nu_n$ for some $n\in \mathbb{N}$ and $\nu\in\{+,-\}$, the $S^\nu_n$-solutions has same humps and the first hump is symmetric around $\frac 1{2n}$, and
     $u$ is a $S^\nu_n$ solution of (1.4) with $\kappa>0$ implies $\lambda\in (0, n^2\pi^2)$.
     Section 3 is devoted to show that (1.4) with $\kappa>0$ has a $S^\nu_n$-solution if and only if
$8n^2 B^2<\lambda<n^2\pi^2$, and complete the proof of Theorem 1.1.
In section 4, we state and prove some properties for the nontrivial solutions of (1.4) with $\kappa<0$. Section 5 is devoted to show that (1.4) with $\kappa<0$ has a $S^\nu_n$-solution if and only if
$n^2\pi^2<\lambda<\infty$, and complete the proof of Theorem 1.2.

\section{The properties of nontrivial solutions of (1.4) with $\kappa>0$}

As the linearization of (1.4) at $(\lambda,u) = (\lambda,0)$ is
given by (1.2), by the local bifurcation theorem of Crandall and
Rabinowitz [2,3], for each $n\in \mathbb{N}$, (1.1) possesses a
curve of $S^\nu_n$-solutions  emanating from  $(\lambda,u) =
(\lambda,0)$ at $\lambda=n^2\pi^2$. Actually, $n^2\pi^2$  is the
unique bifurcation value from  $u =0$ to $S^\nu_n$-solutions  of
(1.4).

\vskip 3mm

The next result shows that any nontrivial solution $u$ of (1.4)
cannot have a degenerate zero.

\vskip 3mm

\noindent{\bf Lemma 2.1} Let $u$ be a nontrivial solution of (1.4)
with $\kappa\neq 0$ for some $\lambda\in (0, \infty)$. Then, $u\in
S^\nu_n$ for some $\nu\in \{+, -\}$ and $n\in \mathbb{N}$.

\vskip 3mm

\noindent{\bf Proof.} Suppose on the contrary that
$$u(\tau)=u'(\tau)=0$$
for some $\tau \in [0, 1]$. Then  $u$ is a
solution of the initial value problem
$$
\left\{
\aligned
&-u''(x)=\lambda u [1+\kappa (u'(x))^2]^{3/2}, \ \
 \ \  \ 0<x<1,\\
&u(\tau)=u'(\tau)=0\\
\endaligned
\right.
\eqno (2.1)
$$
which implies that $u(x)\equiv 0$ for $x\in [0,1]$. This is a contradiction.
\hfill{$\Box$}

\vskip 3mm

The next result shows that this bifurcation is sub-critical.

\vskip 3mm

\noindent{\bf Lemma 2.2} Let us assume $\kappa>0$. Let $u$ be a $S^\nu_n$-solution of (1.1) for some $\lambda\in (0, \infty)$. Then
$0<\lambda < n^2\pi^2$. Moreover,

\noindent(1) All the positive bumps of $u$ have the same shape, and any such bump $B$ satisfies

\ (i) $B$ is symmetric about its mid-point $m_B$;

(ii) $u'$ is strictly decreasing on $B$, so $B$ contains exactly one zero of $u'$
 at $m_B$;

\noindent(2) All the negative bumps of $u$ have the same shape, and any such bump $D$ satisfies

\ (i) $D$ is symmetric about its mid-point $m_D$;

(ii) $u'$ is strictly increasing on $D$, so $m_D$ contains exactly one zero of $u'$
 at $m_D$;

\noindent(3) The positive and negative bumps have the same shape.

  \vskip 2mm

\noindent{\bf Proof.} Without loss of generality, we assume $u\in S^+_n$.
Let
$$b:=u'(0).$$
 Then $b>0$. Notice that $u$ is also a nontrivial solution of  the initial value problem
$$
\left\{
\aligned
&-\Big(\frac{u'(x)}{\sqrt{1+\kappa (u'(x))^2}}\Big)'=\lambda u(x), \ \
 \ \ \ 0<x<1,\\
&u(0)=0, \ \ \  u'(0)=b.\\
\endaligned
\right.
\eqno (2.2)
$$
Denote
$$v=u'.
\eqno (2.3)
$$
Then (1.3) can be written in the form
$$\frac{v'}{(1+\kappa v^2)^{3/2}}=\lambda u
$$
and, hence
$$\frac{vv'}{(1+\kappa v^2)^{3/2}}=\lambda uu'
$$
or, equivalently,
$$\frac d{dx}\Big(-\frac 1{\kappa}(1+\kappa v^2)^{-1/2}+\lambda \frac {u^2}2\Big)=0.
$$
Therefore, for every $x\in [0, 1]$, we have that
$$\lambda \frac {\kappa}2u^2(x)-\frac 1{\sqrt{1+\kappa v^2(x)}}= -\frac 1{\sqrt{1+\kappa b^2}}.
\eqno (2.4)
$$
It now follows from (2.4) that if $b = 0$ then
  $u\equiv 0$ on $[0, 1]$. Hence, if $u\not\equiv 0$ then any zeros of $u$ are simple, and $u\neq 0$ at any zero of $u'$.

  Let $\tau\in [0,1)$ be a zero of $u$.
  Suppose that $u>0$ on a maximal interval $P$ with  $\tau$ is as a left end (similar arguments apply to maximal intervals $N$ on which $u < 0$). Then from the equation in (2.2),
  $$
\frac{u'(x)}{\sqrt{1+\kappa (u'(x))^2}}=\frac{u'(\tau)}{\sqrt{1+\kappa (u'(\tau))^2}}
-\lambda \int ^x_\tau u(s)ds, \ \
 \ \ \ x\in P.\\
\eqno (2.5)
$$
Since the function
 $$\psi(s):=\frac s{\sqrt{1+\kappa s^2}}$$
 is strictly increasing on $(-\infty, \infty)$, it follows from (2.5) that $u'$ is strictly decreasing on $P$. Hence, $P$ contains at most one zero of $u'$.

  Now denote that $P = (x_l, x_r)$, $m := (x_r+x_l)/2$.  Then, by (2.4),
  $$u(x_l)=u(x_r)=0, \ \ \ \ u'(x_l)=-u'(x_r)>0.
  \eqno (2.6)
  $$
  By (2.6), it is easy to check that  both $u(x)$ and $u(x_r+x_l-x)$ are solutions of initial value problem
  $$
\left\{
\aligned
&-\Big(\frac{w'(x)}{\sqrt{1+\kappa (w'(x))^2}}\Big)'=\lambda w(x), \ \
 \ \ \ x_l<x<x_r,\\
&w(x_l)=0, \ \ \  w'(x_l)=u'(x_l).\\
\endaligned
\right.
$$
 Hence, by the above uniqueness result, the solution curves on the intervals $[x_l,m]$, $[m, x_r]$ are symmetric, and $u'(m)=0$. These results show that any positive (and negative) bumps have the properties (i) and (ii) described in the theorem.

  To prove (3), we take $n\geq 2$. Let
  $$0=\tau_0<\tau_1<\cdots<\tau_n=1
  $$
  be the zeros of $u$ in $[0,1]$.
  Let us consider the initial value problem
  $$
\left\{
\aligned
&-v''(x)=\lambda u [1+\kappa (v'(x))^2]^{3/2}, \ \
 \ \  \ 0<x<1,\\
&v(0)=0, \ \ \ v'(0)=u'(0).\\
\endaligned
\right.
\eqno (2.7)
$$
  Since the equation in (1.3) is autonomous, it is easy to check that
   the  function
  $$v(x):=-u(x-\tau_1), \ \ \ \ x\in [\tau_1, \tau_2]
  $$
 is also a solution of the equation in (2.7). Now by the uniqueness results for the
 initial value problem (2.7),
  $u\equiv v$ in $[\tau_0, \tau_1]$, and accordingly  $u\equiv v$ in $[0,1]$, which implies the first positive bump and the first negative bump
  have the same shape.

  \vskip 3mm

  To show that $0<\lambda<n^2\pi^2$, multiplying the equation in (2.1)
  by $\sin(n\pi x)$ and integrating in $(0, \frac 1n)$, we are lead to
  $$
  \aligned
  n^2\pi^2\int^{1/n}_0 \sin (n\pi x) u(x)dx
  &=- \int^{1/n}_0 \sin (n\pi x) u''(x)dx\\
  &=\lambda\int^{1/n}_0 \sin (n\pi x) u(x)V(x)dx\\
  &>\lambda\int^{1/n}_0 \sin (n\pi x) u(x)dx,\\
  \endaligned
  $$
  where $V:=[1+(u')^2]^{3/2}$. Because $V \succ 1$, for as $u = 0$ if $V = 1$. Therefore,
  $\lambda <n^2 \pi^2$. The proof is complete.
  \hfill{$\Box$}

\vskip 3mm

\section{Interval of $\lambda$ in which (1.4) with $\kappa>0$ has $S^\nu_n$-solutions}

Recall that
$$
B:=\int^1_0 \frac{d\theta}{\sqrt{\theta^{-4}-1}}.
$$

\vskip 3mm

\noindent{\bf Lemma 3.1} [1, Theorem 3.2.] The set
$$\mathcal{S}_1:=\{(\lambda, ||u_\lambda||_\infty)\,|\, (\lambda, u_\lambda)\
\text{is a symmetric positive solutions of (1.4)}\}$$
is a
differentiable monotone curve, which connects $(8B^2,\frac
{1}{2B\sqrt{\kappa}})$ with $(\pi^2, 0)$. For each $\lambda\in
(8B^2,\pi^2)$, the problem (1.4) admits a unique symmetric positive
solution. Moreover, if $u_\lambda$ stands for the unique symmetric
positive solution of (1.4), then $u'_\lambda(0)$ is strictly
decreasing in $(8B^2,\pi^2)$ and
$$u'_\lambda(0)=\infty.
$$
Furthermore, the point-wise limit
$$u_{8B^2}:=\underset{\lambda \downarrow 8B^2}\lim \; u_\lambda$$
provides us with a solution of (1.4) at $8B^2$ which is
singular at $0$ and $1$ in the sense that
$$u'_{8B^2}(0)=-u'_{8B^2}(1)=\infty.
$$

 In this section, we shall extend the above result to the case of nodal solutions.

\noindent{\bf Theorem 3.2 } For each $n\in N$ and $\nu\in\{+, -\}$,
the set
$$\mathfrak{S}^\nu_n:=\{(\lambda, ||u_\lambda||_\infty)\,|\, (\lambda, u_\lambda)\
\text{is a} \ S^\nu_n\text{-solutions of}\ (1.4)\}$$ is a
differentiable monotone curve, which connects $(8n^2 B^2,\frac
{1}{2B\sqrt{\kappa}})$ with $(\pi^2n^2, 0)$. For each $\lambda\in
(8n^2 B^2,n^2\pi^2)$, the problem (1.4) admits a
unique ${S}^\nu_n$-solution. Moreover, if $u_\lambda$ stands
for the unique ${S}^\nu_n$-solution of (1.4), then
$u'_\lambda(0)$ is strictly decreasing in $(8n^2 B^2,n^2\pi^2)$ and
$$\underset{\lambda\downarrow 8n^2B^2}\lim\, u'_\lambda(0)=\infty.
$$
Furthermore, the point-wise limit
$$u_{8n^2B^2}:=\underset{\lambda \downarrow 8n^2B^2}\lim \; u_\lambda$$
provides us with a solution of (1.4) at $8n^2B^2$ which is
singular at $0$ and $1$ in the sense that
$$(-1)^j u'_{8n^2B^2}(\frac jn)=\infty, \ \ \ \ \ j\in\{0, 1, \cdots, n\}.
$$

\noindent{\bf Proof.}  By Lemma 2.2, to study the nodal solutions of
(1.1) in $\mathfrak{S}^\nu_n$, it is enough to study the positive
solutions of
$$
\left\{
\aligned
&-\Big(\frac{u'(x)}{\sqrt{1+\kappa (u'(x))^2}}\Big)'=\lambda u(x), \ \
 \ \ \ 0<x<\frac 1n,\\
&u(0)= u(\frac 1n)=0.\\
\endaligned
\right.
\eqno (3.1)
$$
The change of variable
$$u(x)=v(y), \ \ \ \ y=xn, \ \ \  0\leq x\leq \frac 1n
\eqno (3.2)$$
transforms (3.1) into
$$
\left\{
\aligned
&-\Big(\frac{v'(y)}{\sqrt{1+\tilde\kappa (v'(y))^2}}\Big)'=\tilde \lambda v(y), \ \
 \ \ \ 0<y<1,\\
&v(0)= v(1)=0,\\
\endaligned
\right.
\eqno (3.3)
$$
where we are denoting
$$\tilde \kappa:=\kappa\, n^2 , \ \ \ \ \ \
\tilde \lambda:=\frac{\lambda}{n^2}. \eqno (3.4)
$$
As (3.3) is of the same type as (1.4), by the analysis already done
in Section 2, it becomes apparent that (3.1) possesses a positive
symmetric solution around in $\frac1 {2n}$ if and only if
$$8B^2<\tilde\lambda<\pi^2\ \Leftrightarrow\ 8n^2B^2<\lambda<n^2\pi^2.
\eqno (3.5)
$$
Furthermore, we may use Lemma 3.1 to get the desired results.
\hfill{$\Box$}

\vskip 3mm

\section{The properties of nontrivial solutions of (1.4) with $\kappa<0$}

As the linearization of (1.4) at $(\lambda,u) = (\lambda,0)$ is
given by (1.2), by the local bifurcation theorem of Crandall and
Rabinowitz [2,3], for each $n\in \mathbb{N}$ and $\nu\in\{+,-\}$,
(1.4) possesses a curve of $S^\nu_n$-solutions emanating from
$(\lambda,u) = (\lambda,0)$ at $\lambda=n^2\pi^2$. Actually,
$n^2\pi^2$  is the unique bifurcation value from  $u =0$ to
$S^\nu_n$-solutions  of (1.4).

\vskip 3mm

\vskip 3mm

The next result shows that this bifurcation is sup-critical in the
case $\kappa<0$.

\vskip 3mm

\noindent{\bf Lemma 4.1}  Let $\kappa<0$. Let $u$ be a
$S^\nu_n$-solution of (1.4) for some $\lambda\in (0, \infty)$. Then
$n^2\pi^2<\lambda<\infty$. Moreover,

\noindent(1) All the positive bumps of $u$ have the same shape, and any such bump $B$ satisfies

\ (i) $B$ is symmetric about its mid-point $m_B$;

(ii) $u'$ is strictly decreasing on $B$, so $B$ contains exactly one zero of $u'$
 at $m_B$;

\noindent(2) All the negative bumps of $u$ have the same shape, and any such bump $D$ satisfies

\ (i) $D$ is symmetric about its mid-point $m_D$;

(ii) $u'$ is strictly increasing on $D$, so $m_D$ contains exactly one zero of $u'$
 at $m_D$;

\noindent(3) The positive and negative bumps have the same shape.

  \vskip 2mm

\noindent{\bf Proof.} Using the same method to prove Lemma 2.2, with obvious changes, we may get the that the humps of $u$ have the same shape.

  To show that $n^2\pi^2<\lambda<\infty$, multiplying the equation in (2.1) by $\sin(n\pi x)$ and integrating in $(0, \frac 1n)$, we are lead to
  $$
  \aligned
  n^2\pi^2\int^{1/n}_0 \sin (n\pi x) u(x)dx
  &=-\int^{1/n}_0 \sin (n\pi x) u''(x)dx\\
  &=\lambda\int^{1/n}_0 \sin (n\pi x) u(x)V(x)dx\\
  &<\lambda\int^{1/n}_0 \sin (n\pi x) u(x)dx,\\
  \endaligned
  $$
  where $V:=[1+\kappa(u')^2]^{3/2}$. Because $1\succ V$,
  for as $u = 0$ if $V\equiv 1$. Therefore,
  $\lambda >n^2 \pi^2$. The proof is complete.
  \hfill{$\Box$}

\vskip 3mm

\section{Interval of $\lambda$ in which (1.4)
 with $\kappa<0$ has $S^\nu_n$-solutions}

    We only determine the maximum interval of $\lambda$ in which (1.4)
 with $\kappa<0$ has $S^+_n$-solutions since the other case can be treated similarly.

 \vskip 3mm

From Lemma 4.1, $u\in S^+_n$, the humps of $u$ are same, and the
first hump is symmetric around $\frac 1{2n}$,  and
$$u'(x)>0\ \  \text{if} \ x\in [0, \frac 1{2n}), \ \ \ \ u'(\frac 1{2n})=0, \ \ \ \
         u'(x)<0 \   \text{if} \ x\in (\frac 1{2n},\frac 1{n}].
         \eqno (5.1)
         $$
So we only need to study the maximum interval of $\lambda$ in which
(1.4) with $\kappa<0$  has symmetric positive solutions in $(0,
\frac 1n)$.

\vskip 5mm

First, we consider the special case that $u\in S^+_1$ and establish
the following

\vskip 3mm

\noindent{\bf Theorem 5.1} The set
$$\mathfrak{S}_1^+:=\{(\lambda, ||u_\lambda||_\infty)\,|\, (\lambda, u_\lambda)\
\text{is a} \ S^+_1\text{-solutions of}\ (1.4)\}$$ consists of a
differentiable monotone curve in $[\pi^2, \infty)\times [0, \frac 1{2\sqrt{-\kappa}})$. For
each $\lambda\in (\pi^2,\infty)$, the problem (1.4) admits a unique
$S^+_1$-solution. Moreover, if $u_\lambda$ stands for the unique
$S^+_1$-solution of (1.4), then $u'_\lambda(0)$ is strictly
increasing in $(\pi^2,\infty)$.
Furthermore, the point-wise limit
$$u_{\infty}(t):=\underset{\lambda \uparrow \infty}\lim \; u_\lambda(t),
$$
which satisfies
$$u'_\infty(0)=- u'_\infty(1)=\frac 1{\sqrt{-\kappa}}.
$$

\vskip 3mm

To prove Theorem 5.1, we need some preliminaries.

\vskip 3mm
Denote
$$u_0=u(\frac 12), \ \ \ \ \ v_0=u'(0), \ \ \ \ v=u'.
$$
Thanks to (5.1) with $n=1$, (1.5) can be written in the form
$$-\frac{v'}{(1+\kappa v^2)^{3/2}}=\lambda u.
$$
and hence
$$-\frac{vv'}{(1+\kappa v^2)^{3/2}}=\lambda uu',
$$
or, equivalently,
$$\frac d{dx}\Big(-\frac 1{\kappa}(1+\kappa v^2)^{-1/2}+\lambda \frac {u^2}2\Big)=0.
$$
Therefore, for every $x\in[0, 1]$, we have that
$$
\lambda\frac{\kappa}{2}u^2(x)-\frac 1{\sqrt{1+\kappa v^2(x)}}
=-\frac 1{\sqrt{1+\kappa v_0^2}}
=\lambda\frac{\kappa}{2}u_0^2-1,
\eqno (5.2)
$$
$$v(x)=u'(x)=\sqrt{\frac1{-\kappa}}\sqrt{1-\big\{1-\lambda \frac {\kappa} 2(u_0^2-u^2(x))\big\}^{-2}}
\ \ \ \text{for all} \ x\in [0, \frac 1{2}],
\eqno (5.3)
$$
and
$$v(x)=u'(x)=-\sqrt{\frac1{-\kappa}}\sqrt{1-\big\{1-\lambda \frac {\kappa} 2(u_0^2-u^2(x))\big\}^{-2}}
\ \ \ \text{for all} \ x\in [\frac 1{2}, 1].
\eqno (5.4)
$$
Thus, necessarily
$$
\aligned
\frac 12 &=\sqrt{-\kappa}\int^{1/2}_0 \frac{u'(x)}{\sqrt{1-\big\{1-\lambda \frac {\kappa} 2(u_0^2-u^2(x))\big\}^{-2}}}dx\\
         &=\sqrt{-\kappa}\int^{u_0}_0 \frac{1}{\sqrt{1-\big\{1-\lambda \frac {\kappa} 2(u_0^2-u^2)\big\}^{-2}}}du\\
         &=\sqrt{-\kappa}\int^{1}_0 \frac{u_0}{\sqrt{1-\big\{1-\lambda \frac {\kappa} 2u_0^2(1-\theta^2)\big\}^{-2}}}d\theta\\
         \endaligned
         $$
where $\theta:=\frac{u}{u_0}.$

Equivalently,
$$\frac 12= J(\lambda, u_0),
\eqno (5.5)$$
where $J: D(J)\to \mathbb{R}_+:=(0, \infty)$ is the function defined by
$$
J(\lambda,\xi)=\sqrt{-\kappa}\int^{1}_0 \frac{\xi}{\sqrt{1-\big\{1-\lambda \frac {\kappa} 2\xi^2(1-\theta^2)\big\}^{-2}}}d\theta
\eqno (5.6)
$$
for all $(\lambda,\xi)\in D(J)$, where
$$D(J)=\Big\{(\lambda,\xi)\,|\, (\lambda,\xi)\in \mathbb{R}_+\times \mathbb{R}_+\Big\}.
$$
Owing to (5.3), it can be easily seen that the symmetric positive
solutions of (1.4) are in one-to-one correspondence with the pairs
$(\lambda, u_0)\in D(J)$ satisfying (5.5).

The next result provides us with some important monotonicity properties of the function $J(\lambda,\xi)$.

\noindent{\bf Proposition 5.2} \ (i) For every $\xi>0$, the function
$\lambda\mapsto J(\lambda,\xi)$ is decreasing and
$$\lim_{\lambda\downarrow0}J(\lambda,\xi)=\infty.
\eqno(5.7)
$$

\noindent(ii) For every $\lambda>0$, the function $\varphi(\xi):= J(\lambda,\xi)$ satisfies $\varphi'(\xi)>0$, and, hence, it is increasing.

\noindent{\bf Proof.}\ The proof of (i) is straightforward and, hence, we omit its details. The proof of (ii) is more delicate. From (5.6), differentiating with respect to $\xi$ and rearranging terms, we find that
$$
\begin{aligned}
\varphi'(\xi)&=\sqrt{-k}\int_{0}^{1}\{1-[1-\lambda\frac{k}{2}\xi^{2}(1-\theta^{2})]^{-2}\}^{-\frac{3}{2}}[1-\lambda\frac{k}{2}\xi^{2}(1-\theta^{2})]^{-3}\lambda k\xi^{2}(1-\theta^{2})d\theta\\
&+\sqrt{-k}\int_{0}^{1}\{1-[1-\lambda\frac{k}{2}\xi^{2}(1-\theta^{2})]^{-2}\}^{-\frac{1}{2}}d\theta\\
&=\sqrt{-k}\int_{0}^{1}\{1-[1-\lambda\frac{k}{2}\xi^{2}(1-\theta^{2})]^{-2}\}^{-\frac{3}{2}}\{[1-\lambda\frac{k}{2}\xi^{2}(1-\theta^{2})]^{-3}\lambda k\xi^{2}(1-\theta^{2})\\
&\ \ \ \ +1-[1-\lambda\frac{k}{2}\xi^{2}(1-\theta^{2})]^{-2}\}d\theta.
\end{aligned}
$$
Notice that $k<0$
and set
$$\alpha(\theta):=[1-\lambda\frac{k}{2}\xi^{2}
(1-\theta^{2})]^{-1}, ~~ 0\leq\theta\leq1.  \eqno(5.8)$$
Then
$$0<\alpha(\theta)\leq1~~~  \text{and} ~~~ 1-\alpha^{2}(\theta)\geq0~~~ \forall \theta\in[0,1]$$
and
$$\varphi'(\xi)=\sqrt{-k}\int_{0}^{1}(1-\alpha^{2})^{\frac{3}{2}}[\alpha^{3}\lambda k\xi^{2}(1-\theta^{2})+(1-\alpha^{2})]d\theta.$$
On the other hand, according to (5.8), it becomes apparent that
$$\alpha\lambda k\xi^{2}(1-\theta^{2})=2(\alpha-1)$$
and, hence,
$$
\begin{aligned}
\varphi'(\xi)&=\sqrt{-k}\int_{0}^{1}(1-\alpha^{2})^{\frac{3}{2}}[2\alpha^{2}(\alpha-1)+(1-\alpha^{2})]d\theta\\
&=2\sqrt{-k}\int_{0}^{1}(1-\alpha^{2})^{\frac{3}{2}}(1-\alpha)^{2}(\alpha+\frac{1}{2})d\theta>0.
\end{aligned}
$$
The proof is complete. \hfill$\Box$

\vskip 3mm

\noindent{\bf Proof of Theorem 5.1.}  By the fact
$$-\kappa v_0^2<1.$$
This together with the symmetric property of $u$ imply that
$$\xi<\frac 1{2\sqrt{-\kappa}}.
$$

From Proposition 5.2, there exists an interval $(0, \rho)\subseteq
\big(0, \frac 1{2\sqrt{-\kappa}}\big)$, such that for every $\xi\in
(0, \rho)$, there exists a unique $\lambda(\xi)>0$ satisfying
$$J\big(\lambda(\xi),\xi\big)=\frac 12,
$$
i.e.
$$\frac 12=\frac 1{\sqrt{\lambda(\xi)}}\int^1_0
\frac{1-\lambda\frac \kappa 2 \xi^2(1-\theta^2)}
  {\sqrt{1-\theta^2-\lambda \frac \kappa 4\xi^2(1-\theta^2)^2}}d\theta$$
and, therefore, since $\pi^2<\lambda(\xi)$, letting $\xi\to 0$, we get
$$\frac 12=\frac 1{\sqrt{\lambda(0)}}\int^1_0
\frac{d\theta}
  {\sqrt{1-\theta^2}}
  =\frac 1{\sqrt{\lambda(0)}}\frac {\pi}2.
 $$
   Hence, $\lambda(0):=\lim\limits_{\xi\rightarrow 0}\lambda(\xi)=\pi^2.$

\vskip2mm

By the monotonicity and continuity of $J$, $\xi\to \lambda(\xi)$
must be continuous and strictly increasing. The differentiability of
$\mathfrak{S}^+_1$ is a direct consequence from Proposition 5.2(ii)
and the implicit function theorem.

     The function $u'_\lambda(0)$ is strictly increasing
     for $\lambda\in (\pi^2, \infty)$ is an immediate consequence of
the fact $\xi\to \lambda(\xi)$ is strictly increasing in $(0, \rho)$
 and the relation
$$
\lambda\frac{\kappa}{2}u^2(x)-\frac 1{\sqrt{1+\kappa v^2(x)}}
=-\frac 1{\sqrt{1+\kappa v_0^2}}
=\lambda\frac{\kappa}{2}\xi^2-1.
$$

On the other hand, let
$$\underset{\lambda\uparrow \infty}\lim \; u_\lambda(\frac 12)
=u_\infty(\frac 12)=:\xi^*.
$$
Then $\xi^*>0$. Combining this with the fact $-\frac
1{\sqrt{1+\kappa v_0^2}} =\lambda\frac{\kappa}{2}\xi^2-1 $ and
letting $\lambda\to \infty$, we get
$$\underset{\lambda\uparrow \infty}\lim u'_\lambda(0)=u'_\infty(0)=\frac1{\sqrt{-\kappa}},
$$
and this completes the proof. \hfill{$\Box$}

\vskip 3mm

Now, we are in the position to prove Theorem 1.2.

\vskip 3mm

\noindent{\bf {Proof of Theorem 1.2.}}  By Lemma 4.1, to study the
set of  $S^\nu_n$-solutions of (1.4), it is enough to
study the set of  positive solutions of
$$
\left\{
\aligned
&-\Big(\frac{u'(x)}{\sqrt{1+\kappa (u'(x))^2}}\Big)'=\lambda u(x), \ \
 \ \ \ 0<x<\frac 1n,\\
&u(0)= u(\frac 1n)=0.\\
\endaligned
\right. \eqno (5.9)
$$
The change of variable
$$u(x)=v(y), \ \ \ \ y=xn, \ \ \  0\leq x\leq \frac 1n,
\eqno (5.10)$$ transforms (5.9) into
$$
\left\{
\aligned
&-\Big(\frac{v'(y)}{\sqrt{1+\tilde\kappa (v'(y))^2}}\Big)'=\tilde \lambda v(y), \ \
 \ \ \ 0<y<1,\\
&v(0)= v(1)=0,\\
\endaligned
\right. \eqno (5.11)
$$
where we are denoting
$$\tilde \kappa:=\kappa\, n^2 , \ \ \ \ \ \ \tilde \lambda:=\lambda/n^2.
\eqno (5.12)
$$
Notice that all of the conclusions of Theorem 5.1 for (5.11) are
still valid for (5.9) via the transformation (5.10). Thus, Theorem 1.2 can be deduced from (5.12)
and Theorem 5.1.
 \hfill{$\Box$}

\vskip 10mm

\centerline {\bf REFERENCES}\vskip5mm\baselineskip 18pt
\begin{description}

\item{[1]}  S. Cano-Casanova, J. L\'{o}pez-G\'{o}mez, K. Takimoto, A quasilinear parabolic perturbation of the linear heat equation, J. Differential Equations 252 (2012)
323-343.

\item{[2]}   M. G. Crandall, P. H. Rabinowitz, Bifurcation from simple eigenvalues, J. Funct. Anal. 8 (1971) 321-340.

\item{[3]} M. G. Crandall, P. H. Rabinowitz, Bifurcation, perturbation from simple eigenvalues and linearized stability, Arch. Ration. Mech. Anal. 52 (1973) 161-180.

\item{[4]}   J. L\'{o}pez-G\'{o}mez,Spectral Theory and Nonlinear Functional Analysis,  Chapman \& Hall/CRC Res. Notes Math., vol. 426, Chapman \& Hall/CRC, Boca Raton, FL, 2001.

\item{[5]}   J. L\'{o}pez-G\'{o}mez, C. Mora-Corral, Algebraic Multiplicity of Eigenvalues of Linear Operators, Oper. Theory Adv. Appl., vol. 177,
    Birkh\"{a}user/Springer, Basel, Boston/Berlin, 2007.

\item{[6]}  M. Nakao, A bifurcation problem for a quasi-linear elliptic boundary value problem, Nonlinear Anal. TMA 14 (1990) 251-262.

\item{[7]}\ S. Cano-Casanova, J.  L\'{o}pez-G\'{o}mez, K. Takimoto,  A weighted quasilinear equation related to the mean curvature operator, Nonlinear Anal. 75(15) (2012) 5905-5923.

\item{[8]} D.Gilbarg, N.S.Trudinger, Elliptic Partial Differential Equations of Second Order,\\  Springer-Verlag, Berlin, 2001.

\item{[9]}\ S.-Y. Cheng, S.-T. Yau, Maximal spacelike hypersurfaces in the Lorentz-Minkowski spaces, Ann. of Math. 104 (1976) 407-419.

\item{[10]}\ R. Bartnik, L. Simon, Spacelike hypersurfaces with prescribed boundary values and mean curvature, Comm. Math.
Phys. 87 (1982-1983) 131-152.

\item{[11]}\  C. Bereanu, P. Jebelean, J. Mawhin, Radial solutions for Neumann problems involving mean curvature operators in
Euclidean and Minkowski spaces, Math. Nachr. 283 (2010) 379-391.

\item{[12]}\ C. Bereanu, P. Jebelean, P. J. Torres, Positive radial
 solutions for  Dirichlet problems with mean curvature operators in
 Minkowski space, J. Funct. Anal. 264 (2013) 270-287.

\item{[13]}\ C. Bereanu, P. Jebelean, P. J. Torres, Multiple positive
 radial solutions for a Dirichlet problem involving the mean curvature
 operator in Minkowski space, J. Funct. Anal. 265(4) (2013) 644-659.

\item{[14]}\ M. F. Bidaut-V\'{e}ron, A. Ratto, Spacelike graphs with prescribed mean curvature, Differential Integral Equations 10 (1997), 1003-1017.

\item{[15]}\ I. Coelho, C. Corsato, F. Obersnel, P. Omari,  Positive solutions of the Dirichlet problem for the one-dimensional Minkowski-curvature equation, Adv. Nonlinear Stud. 12(3) (2012)  621-638.

\item{[16]}\ R. L\'{o}pez, Stationary surfaces in Lorentz-Minkowski space, Proc. Roy. Soc. Edinburgh Sect. A 138A (2008) 1067-1096.

\item{[17]}\ J. Mawhin, Radial solution of Neumann problem for periodic perturbations of the mean extrinsic curvature operator,
Milan J. Math. 79 (2011) 95-112.

\item{[18]}\ A. E. Treibergs, Entire spacelike hypersurfaces of constant mean curvature in Minkowski space, Invent. Math. 66 (1982) 39-56.

\item{[19]}\ A. \ Azzollini, \ Ground state solution for a problem with mean curvature operator in Minkowski space, J. Funct. Anal. 266 (2014), 2086-2095.

\item{[20]}\  H. Pan, R. Xing, Sub- and supersolution methods for prescribed mean curvature equations with Dirichlet boundary conditions. J. Differential Equations 254 (2013), 1464-1499.

\end{description}
\end{document}